\documentclass{amsart}
\usepackage{amsrefs}
\usepackage{color}
\newtheorem{theorem}{\bf Theorem}[section]

\newtheorem{lemma}{\bf Lemma}[section]
\newtheorem{definitions}{\bf Definition}[section]
\newtheorem{proposition}{\bf Proposition}[section]

\newtheorem{definition}{\bf Definition}[section]

\newcommand{\um}{\ua}
\newcommand{\ua}{u_\alpha}
\newcommand{\tva}{\tilde{v}_\alpha}
\newcommand{\xm}{x_\alpha}

\newcommand{\R}{{\mathbb R}}
\newcommand{\rn}{{\mathbb R}^n}
\newcommand{\rnm}{{\mathbb R}^n_-}
\newcommand{\N}{{\mathbb N}}
\newcommand {\crit} {2_k^\sharp}

\begin{document}

\title [Struwe's Decomposition for a Polyharmonic Operator]{Struwe's Decomposition  for a Polyharmonic Operator on a Compact 
Riemannian Manifold with or without boundary}
\author{Saikat Mazumdar}
\address{Institut Elie Cartan de Lorraine, Universit\'e de Lorraine, BP 70239, 54506    Vand\oe uvre-l\`es-Nancy, France}
\email{saikat.mazumdar@univ-lorraine.fr}

\begin{abstract} Given a high-order elliptic operator on a compact manifold with or without boundary, we perform the decomposition of Palais-Smale sequences for a nonlinear problem as a sum of bubbles. This is a generalization of the celebrated 1984 result of Struwe \cite{struwe.1984}. Unlike the case of second-order operators, bubbles close to the boundary might appear. Our result includes the case of a smooth bounded domain of $\rn$.
\end{abstract}
\date{March 25th, 2016}
\thanks{This work is part of the PhD thesis of the author, funded by "F\'ed\'eration Charles Hermite" (FR3198 du CNRS) and "R\'egion Lorraine". The author acknowledges these two institutions for their supports.}

\maketitle

\section{Introduction}
Let $({M},g)$ be a smooth, compact Riemannian manifold of dimension $n$  with or without boundary. In the latter case  
 we understand that $\overline{M}$ is a compact,   oriented submanifold of $ (\tilde{M},g)$ which is itself a  smooth, compact Riemannian manifold  without boundary and with the same metric $g$. As one checks, this includes smooth bounded domains of $\rn$. When the boundary $\partial M\neq \emptyset$, we let $\nu$ be its outward oriented normal vector in $\tilde{M}$. Let $k$ be a positive integer such that $2k < {n}$. We define the Sobolev space $H^{2}_{k,0}(M)$ as the completion of $C^\infty_c(M)$ for the norm $u\mapsto \sum_{i=0}^k\Vert\nabla^{i}u\Vert_2$. This norm is equivalent (see Robert \cite{robert.gjms}) to the Hilbert norm $\Vert u\Vert_{H_{k}^2} := \left(
 \sum_{l =0}^{k} \int_M   (\Delta_g^{{l}/2} u )^2 ~{dv_g}\right)^{1/2}$ where $\Delta_g:=-\hbox{div}_g(\nabla)$ is the Laplace-Beltrami operator and, for $\alpha$ odd, $\Delta_{g}^{\alpha} u \Delta_g^\alpha v:=(\nabla \Delta_{g}^{\frac{\alpha-1}{2}} u,\nabla \Delta_{g}^{\frac{\alpha-1}{2}} v)_g $ for all $u,v\in H_{k}^2(M)$. For details we refer to Aubin \cite{aubin.book} and Hebey \cite{hebey.book}.

\medskip\noindent We consider the functional
 \begin{equation*}
 I(u):= \frac{1}{2} \int_M (\Delta^{k/2}_g u)^2~ dv_{g}+\frac{1}{2}\sum_{l=0}^{k-1}\int_M A_{l}(\nabla^lu,\nabla^l u)\, dv_g
 - \frac{1}{2_k^{\sharp}} \int_M| u |^{ 2_k^{\sharp}}~{dv_g} 
 \end{equation*}
where for all $l \in \{ 0, \ldots, k-1\}$, $A_{l} $  is a smooth $T_{2l}^{0}$-tensor field on $M$ and $A_{l} $ is symmetric (that is $A_{l}(X,Y)  =A_{l}(Y,X)  $ for all $T_{0}^{l}$-tensors $X,Y$ on $M$).  Here, $\crit:=\frac{2n}{n-2k}$ is the critical Sobolev exponent such that $H_{k,0}^2(M)\hookrightarrow L^{\crit}(M)$ is continuous, which makes the definition of $I$ consistent for all $u\in H_{k,0}^2(M)$. Critical points $u\in H_{k,0}^2(M)$ for $I$ are weak solutions to the pde
\begin{eqnarray}{\label{eqn1}}
 \left \{ \begin{array} {ll}
          Pu  = |u|^{2_{k}^{\sharp}-2} u & \text{in } ~ {M}\\
          \partial^{\alpha}_{\nu} u=0  &\text{on } ~\partial {M} \quad \text{for} \ \ |\alpha| \leq k-1            \end{array} \right. 
\end{eqnarray} 
where for any $u\in C^{2k}(M)$, we define
\begin{align*}
 Pu:= \Delta_{g}^{k}u + \sum_{l=0}^{k-1} (-1)^{l}\nabla^{j_{l} \ldots j_{1}} \left( (A_{l})_{i_{1} \ldots i_{l}, j_{1} \ldots j_{l}} \nabla^{i_{1} \ldots i_{l}} u\right)
 \end{align*}
and where we say that $u \in H^{2}_{k,0}(M)$ is  a weak solution to \eqref{eqn1} if  
\begin{align*}
  \int_M\Delta_g^{k/2} u , \Delta_g^{k/2} \varphi  ~{dv_g} + 
  \sum_{l =0}^{k-1} \int_M
   A_{l} (\nabla^{l} u\nabla^{l} \varphi) ~dv_{g}
     = \int_M \left| u \right|^{ 2_k^{\sharp}-2}u \varphi ~{dv_g} 
\end{align*}
for all $\varphi \in H_{k,0}^2(M)$. As shown by the regularity theorem in Mazumdar \cite{mazumdar.gjms}, a weak solution $u$ to \eqref{eqn1} is indeed a  strong solution,  $u\in C^{2k}(\overline{M})$.

\begin{definitions} Let $(X,\Vert\cdot\Vert)$ be a Banach space and $F\in C^1(X)$. A sequence $(u_\alpha)$ in $X$ is  said to be a Palais-Smale sequence for $F$ if $(F(u_\alpha))_\alpha$ has a limit in $\R$ when $\alpha\to +\infty$, while $  DF(u_\alpha) \rightarrow 0$ strongly in $X^\prime$ as $\alpha \rightarrow + \infty$. 
\end{definitions}
\noindent In this paper, we describe the lack of relative compactness of Palais-Smale sequences for $I$, which is due to the noncompact embedding $H_{k,0}^2(M)\hookrightarrow L^{\crit}(M)$. For $\Omega$ any open domain of $\rn$, we let  ${\mathcal D}_k^2(\Omega)$ be the completion of $C^\infty_c(\Omega)$ for the norm $u\mapsto\Vert\Delta^{k/2}u\Vert_2$. The limiting equations of \eqref{eqn1} are
\begin{equation}\label{lim:1}
 \Delta^{k} u = \left| u\right|^{2_{k}^{\sharp}-2} u ~ \text{in} ~\R^{n},\; u\in {\mathcal D}_k^2(\rn)
\end{equation} 
\begin{equation}\label{lim:2}
\left\{\begin{array}{ll}
\Delta^{k} u = \left| u\right|^{2_{k}^{\sharp}-2} u &\hbox{in }\rnm\\
\partial^{\alpha}_{\nu} u=0  &\hbox{on }\partial\rnm
\end{array}\right\},\; u\in {\mathcal D}_k^2(\rnm)
\end{equation} 
where $\Delta:=\Delta_{\hbox{Eucl}}$ is the Laplacian on $\rn$ endowed with the Euclidean metric Eucl. Associated to the functional $I$ is the limiting functional
\begin{equation*}\label{def:E}
E(u):=\frac{1}{2}\int_{\rn}(\Delta u)^2\, dx-\frac{1}{\crit}\int_{\rn}|u|^{\crit}\, dx\hbox{ for all }u\in {\mathcal D}_k^2(\rn).
\end{equation*}
Our main theorem below shows that the lack of convergence to a solution of  equation \eqref{eqn1} is described by a sum of Bubbles:

\begin{theorem}\label{th:1}
Let $(\ua)$ be  a Palais-Smale sequence for the functional $I$ on the space $ H_{k,0}^2(M)$. Then there  exists  $d\in\N$ bubbles $[(x_\alpha^{(j)}), (r_\alpha^{(j)}), u^{(j)}]$, $j=1,...,d$, $($see Definition \ref{def:bubbles} below$)$  there exists $u_\infty \in H_{k,0}^2(M) $ a solution to $\eqref{eqn1}$ such that, up to a subsequence,
\begin{equation*}
\ua =u_\infty + \sum_{j=1}^d B_{x_\alpha^{(j)},r_\alpha^{(j)}}(u^{(j)})+o(1)\hbox{ where }\lim_{\alpha\to +\infty}o(1)=0\hbox{ in }H_{k,0}^2(M)
\end{equation*}
and
\begin{equation*}
\displaystyle  I(\ua)= I(u_\infty) +  \sum_{j=1}^d ~E(u^{(j)}) + o(1)  \quad \text{as $\alpha \rightarrow + \infty$.}
\end{equation*}
\end{theorem}
\noindent
In Section \ref{sec:bubbles}, Bubbles are defined up to a term going to $0$ strongly, which is relevent here. As one checks, given $u\in {\mathcal D}_k^2(\rn)$ a nontrivial weak solution to \eqref{lim:1} or \eqref{lim:2}, then multiplying the equation by $u$ and integrating by parts yields
\begin{equation}\label{def:beta}
E(u) \geq \beta^{\sharp} := \frac{k}{n} K_0(n,k)^{-n/2k}
\end{equation}
where $K_0(n,k)$ be the best constant of the embedding ${\mathcal D}_k^{2}(\R^n) \hookrightarrow L^{2_k^{\sharp}}(\R^n)$, that is
\begin{equation}{\label{sinq}}
K_0(n,k)^{-1} = \inf_{u \in {\mathcal D}_k^2(\rn)\backslash \{0\}} \frac{\int_{\rn}(\Delta^{k/2} u)^2\, dx}{\left(\int_{\rn}|u|^{\crit}\, dx\right)^{\frac{2}{\crit}}}
\end{equation}

\medskip\noindent When the Palais-Smale sequence is nonnegative, the bubbles are positive and correspond to positive solutions to \eqref{lim:1}. As shown in Lions \cite{PLL}, Swanson \cite{swanson}, Ge-Wei-Zhou \cite{poly}, these solutions are exactly the extremals for \eqref{sinq} and are  of the form 
\begin{align}\label{eq:ext}
u=U_{a,\lambda}:=  \alpha_{n,k}\left( \frac{\lambda}{1+ \lambda^{2}|\cdot-a|^{2}} \right)^{\frac{n-2k}{2}}~ a \in \R^{n}, \lambda >0
\end{align}
where $\alpha_{n,k}>0$ is explicit. We then get the following:
\begin{theorem}\label{th:2}
Let $(\ua)$ be  a Palais-Smale sequence for the functional $I$ on the space $ H_{k,0}^2(M)$. We assume that $\ua\geq 0$ for all $\alpha\in\N$. Then there exists $u_\infty \in H_{k,0}^2(M) $ a solution to \eqref{eqn1}, there exists $d\in\N$, there exist $(x_\alpha^{(1)}),\dots,(x_\alpha^{(d)})\in M$, $(r_\alpha^{(1)}),\dots,(r_\alpha^{(d)})\in (0,+\infty)$ such that $r_\alpha^{(j)}\to 0$ and $r_\alpha^{(j)}=o(d(x_\alpha^{(j)},\partial M))$ as $\alpha\to +\infty$ for all $j=1,...,d$, such that, up to a subsequence,
\begin{equation*}
\ua =u_\infty + \sum_{j=1}^d \eta\left((\tilde{r}_\alpha^{(j)})^{-1}\hbox{exp}_{x_\alpha^{(j)}}^{-1}(\cdot)\right)\alpha_{n,k}\left(\frac{r_\alpha^{(j)}}{(r_\alpha^{(j)})^2+d_g(\cdot,x_\alpha^{(j)})^2}\right)^{\frac{n-2k}{2}}+o(1)
\end{equation*}
where $\lim_{\alpha\to +\infty}o(1)=0\hbox{ in }H_{k,0}^2(M)$, and $\eta$ and  $(\tilde{r}_\alpha^{(j)})'s$ are as in \eqref{def:rt}. Moreover, 
\begin{equation*}
\displaystyle  I(\ua)= I(u_\infty) +  d \beta^\sharp + o(1)  \quad \text{as $\alpha \rightarrow + \infty$}
\end{equation*}
where $\beta^\sharp$ is as in \eqref{def:beta}.
\end{theorem}

\medskip\noindent When $k=1$ and $M$ is a smooth bounded domain of $\rn$, Theorem \ref{th:1} is the pioneering  result of Struwe \cite{struwe.1984}. There have been several extensions. Without being exhaustive, we refer to Hebey-Robert \cite{Hebey-robert} for $k=2$ and manifolds without boundary, Saintier \cite{saintier} for the $p-$Laplace operator, El-Hamidi-V\'etois \cite{vetois.arma} for anisotropic operators and Almaraz \cite{almaraz} for nonlinear boundary conditions. A general reference for description as bubbles is the monograph by Fieseler-Tintarev \cite{tintarev.book}. Another possible description is in the sense of measures as in Lions \cite{PLL}: a general result of this flavour for high order elliptic operators on manifolds is in Mazumdar \cite{mazumdar.gjms}.

\smallskip\noindent Palais-Smale sequence are produced via critical point techniques, like the Mountain-Pass Lemma of Ambrosetti-Rabinowitz \cite{a.r} or other topological methods (see for instance the monograph Ghoussoub \cite{ghoussoub.book} and the references therein). Concerning higher-order problems, we refer to Bartsch-Weth-Willem \cite{bww}, Ge-Wei-Zhou \cite{poly}, Mazumdar \cite{mazumdar.gjms}, the general monograph Gazzola-Grunau-Sweers \cite{poly.book} and the references therein. Theorem \ref{th:1} is used by the author in \cite{mazumdar.gjms} to get Coron-type solutions to equation \eqref{eqn1}.

\subsection*{Acknowledgements}

I would like to express my deep gratitude to Professor Fr\'ed\'eric Robert  and Professor Dong Ye, my thesis  supervisors, for their patient guidance, enthusiastic encouragement and useful critiques of this work.

\section{Definition of Bubbles}\label{sec:bubbles}
\noindent
In the spirit of the exponential map, we first cook up a chart around any boundary point. We fix $x_0\in\partial M$. Since $M$ is a smooth submanifold of $\tilde{M}$, there exist $\Omega$ an open subset of $\tilde{M}$ with $x_0\in \Omega$, there exists $U\subset\rn$ open with $0\in U$, such that for any $x\in \Omega\cap\partial M$ there exists $\mathcal{T}_{x}\in C^\infty(U,\tilde{M})$ having the following properties.

\begin{equation}\label{def:T}
\left\{\begin{array}{ll}
\bullet & \mathcal{T}_{x}(0)=x\\
\bullet & \mathcal{T}_x\hbox{ is a smooth diffeomorphism onto its image }\mathcal{T}_x(U).\\
\bullet & \mathcal{T}_{x} \left( U \cap \{x_{1} <0 \} \right)= \mathcal{T}_{x} (U) \cap M\\
\bullet & \mathcal{T}_{x} \left( U \cap \{x_{1} =0 \} \right)= \mathcal{T}_{x} (U) \cap \partial M\\
\bullet & (x,z)\mapsto \mathcal{T}_x(z)\hbox{ is smooth from }\Omega\times U\hbox{ to }\tilde{M}\\
\bullet & d{\mathcal T}_x(0): \rn\to T_{x} M\hbox{ is  an isometry}\\
\bullet & d{\mathcal T}_x(0)[e_1]=\nu_{x}\hbox{ where }\nu_{x}\hbox{ is the outer unit normal vector to }\partial M\\
 & \hbox{ at the point }x. 
\end{array}\right.
\end{equation}
This map is defined uniformly with respect to $x$ in a neighborhood $\Omega$ of a fixed point $x_0\in\partial M$. By a standard abuse of notation, we will always consider $x\mapsto {\mathcal T}_x$ without any reference to $\Omega$ or $x_0$: this will make sense in the sequel since the relevant points will always be in the neighborhood of a fixed point.

\begin{definition}\label{def:bubbles} A ``Bubble'' is a triplet $\left[(x_\alpha), (r_\alpha), u\right]$ where $x_\alpha\in \overline{M}$ is a convergent sequence,  $r_\alpha>0$ for all $m\in\N$ with $\lim_{\alpha\to +\infty}r_\alpha=0$ and 
\begin{align*}
&\hbox{either }\left\{x_\alpha\in M,\; \lim_{\alpha\to +\infty}\frac{d(x_\alpha,\partial M)}{r_\alpha}=+\infty\hbox{ and }u\in {\mathcal D}_k^2(\rn)\hbox{ satisfies \eqref{lim:1}}\right\}\\
&\hbox{ or }\left\{x_\alpha\in \partial M\hbox{ and }u\in {\mathcal D}_k^2(\rnm)\hbox{ satisfies \eqref{lim:2}}\right\}
\end{align*}
If $x_\alpha\in M$, we let $\tilde{r}_\alpha>0$ be such that 
\begin{equation}\label{def:rt}
\lim_{\alpha\to +\infty}\tilde{r}_\alpha=\tilde{r}_\infty\in \left[0,\frac{i_g(\tilde{M})}{2}\right)\; , \; \lim_{\alpha\to +\infty}\frac{r_\alpha}{\tilde{r}_\alpha}=0\hbox{ and }\tilde{r}_\alpha<\frac{d_g(x_\alpha,\partial M)}{2}
\end{equation}
and we define
$$B_{x_\alpha,r_\alpha}(u):=\eta\left(\frac{\hbox{exp}_{x_\alpha}^{-1}(x)}{\tilde{r}_\alpha}\right)r_\alpha^{-\frac{n-2k}{2}}u\left(\frac{\hbox{exp}_{x_\alpha}^{-1}(x)}{r_\alpha}\right)$$
where $\eta\in C^\infty_c(B_0(i_g(\tilde{M})))$ is identically $1$ in a neighborhood of $0$. Here, the exponential map is taken on the ambient manifold $(\tilde{M},g)$.

\smallskip\noindent If $x_\alpha\in\partial M$, we let $x_0:=\lim_{\alpha\to +\infty}x_\alpha$, and we define
$$B_{x_\alpha,r_\alpha}(u):=\eta\left(\mathcal{T}_{x_\alpha}^{-1}(x)\right)r_\alpha^{-\frac{n-2k}{2}}u\left(\frac{\mathcal{T}_{x_\alpha}^{-1}(x)}{r_\alpha}\right)$$
where $\mathcal{T}_{x}$ is as in \eqref{def:T}, $\Omega$ is a neighborhood of $x_0\in\partial M$ and $\eta\in C^\infty_c(U)$ is identically $1$ in a neighborhood of $0$. 
\end{definition}

\medskip\noindent Beside $\left[(x_\alpha), (r_\alpha), u\right]$, the definition of a bubble depends on the choice of the cut-off function $\eta$, the radius $\tilde{r}_\alpha$ and the chart ${\mathcal T}_x$. However, as shown in the proposition below, after quotienting by sequences going to $0$, the class of a Bubble is independent of these later parameters.

\begin{proposition}\label{prop:uniqueness} The definition of Bubbles depend only on $\left[(x_\alpha), (r_\alpha), u\right]$, up to a sequence going to $0$ strongly in $H_{k,0}^2(M)$.
\end{proposition}

\smallskip\noindent{\it Proof of Proposition \ref{prop:uniqueness}.} We first assume that $u\in {\mathcal D}_k^2(\rn)$  satisfies \eqref{lim:1} and that
\begin{equation}\label{case:1}
\lim_{\alpha \to + \infty}\frac{d_{g}(x_{\alpha}, \partial M)}{r_{\alpha}}= + \infty.
\end{equation}
For $i=1,2$, we set the bubbles $B^{i}_{\alpha}:=\eta^{i}\left((\tilde{r}^{i}_\alpha)^{-1}\hbox{exp}_{x_\alpha}^{-1}(\cdot)\right)r_\alpha^{-\frac{n-2k}{2}} u\left(r_\alpha^{-1}\hbox{exp}_{x_\alpha}^{-1}(\cdot)\right)$, where $\eta^{i}\in C^\infty_c(B_0(2a_i))$, $\eta^{i} \equiv 1 $ in $ B_{0}(a_i)$ with $0< 2a_i\leq \iota_{g}(\tilde{M})$;  $\tilde{r}^{i}_{\alpha}>0$ are as in \eqref{def:rt}.  We let $r_{\alpha}^{max}= \max\{ a_1 \tilde{r}_{\alpha}^{1}, a_2 \tilde{r}_{\alpha}^{2}\}$ and  $r_{\alpha}^{min}= \min \{ a_1 \tilde{r}_{\alpha}^{1}, a_2 \tilde{r}_{\alpha}^{2}\}$, and  let $\epsilon^{max}_{\alpha}= r_{\alpha}/{r}_{\alpha}^{max}$ and  $\epsilon^{min}_{\alpha}= r_{\alpha}/{r}_{\alpha}^{min}$. Then $\lim_{ \alpha \to 0} \epsilon_{\alpha}^{max}=0$ and $\lim_{ \alpha \to 0} \epsilon_{\alpha}^{min}=0$. The comparison lemma 9.1 of \cite{mazumdar.gjms} yields  $C>0$ such that for any $R>0$ and $\alpha$ large 
\begin{eqnarray*}
\sum_{l=0}^{k} \Vert \Delta_{g}^{l/2}\left( B^{1}_\alpha-B^{2}_\alpha \right)\Vert_2^2& \leq&    \sum_{l=0}^{k} \int _{B_{ 2 {r}^{max}_{\alpha} } (x_{\alpha}) \backslash  B_{ {r}^{min}_{\alpha} } (x_{\alpha}) }  \left(\Delta_{g}^{l/2}\left( B^{1}_\alpha-B^{2}_\alpha \right) \right)^{2}~ dv_{g}\\
&&\leq  \sum_{i=1,2}\sum_{l=0}^{k} \int_{M\setminus B_{R r_\alpha}(x_\alpha)}  \left(\Delta_{g}^{l/2}\left( B^{i}_{x_\alpha,r_\alpha}(u) \right) \right)^{2}~ dv_{g}.
\end{eqnarray*}
Therefore, using \eqref{ineq:11}, we get that $B^{1}_\alpha-B^{2}_\alpha=o(1)$ in $H_k^2(M)$ as $\alpha\to +\infty$.

\medskip\noindent Now  we consider the case of a  boundary bubble, that is $x_{\alpha} \in \partial M$ and and $u\in {\mathcal D}_k^2(\R^{n}_{-})$ satisfies \eqref{lim:2}. For $i=1,2$, we set $B^{i}_\alpha:=\eta^{i}\left({\mathcal{T}^{1}}_{x_\alpha}^{-1}(\cdot)\right)r_\alpha^{-\frac{n-2k}{2}}u\left(r_\alpha^{-1}{\mathcal{T}^{i}}_{x_\alpha}^{-1}(\cdot)\right)$
where ${\mathcal{T}}"_{x}$, $i=1,2$, are as in \eqref{def:T}, $U$ is a neighborhood of $x_0\in\partial M$ and $\eta^{1}, \eta^{2}\in C^\infty_c(U)$ are identically $1$ in a neighborhood of $0$. One has   
\begin{align*}
& \sum_{l=0}^{k} \int_M \left(\Delta_{g}^{l/2}\left( B^1_\alpha-B^2_\alpha \right) \right)^{2}~ dv_{g} \leq \notag\\
& \sum_{l=0}^{k} \int_{ D_\alpha(R)  \cap M}  \left(\Delta_{g}^{l/2}\left( B^1_\alpha-B^2_\alpha \right) \right)^{2}~ dv_{g} +  \sum_{l=0}^{k} \int_{M \setminus  D_\alpha(R)  }  \left(\Delta_{g}^{l/2}\left( B^1_\alpha-B^2_\alpha \right) \right)^{2}~ dv_{g}
\end{align*}
where $D_\alpha(R):={\mathcal{T}}^{1}_{x_{\alpha}}(B_{0}(r_{\alpha}R) )  \cup   {\mathcal{T}}^{2}_{x_{\alpha}} (B_{0}(r_{\alpha}R))   $
\noindent
It follows as in  the comparison Lemma 9.1  of \cite{mazumdar.gjms}  that there exists  $C>0$ such that for $\alpha$ large
\begin{align*}
& \sum_{l=0}^{k} \int_{ D_\alpha(R)\cap M}  \left(\Delta_{g}^{l/2}\left( B^1_\alpha-B^2_\alpha \right) \right)^{2}~ dv_{g} \leq  \notag \\
&  C\sum_{l=0}^{k} \int_{ \left( B_{0}(r_{\alpha}R) \cup \Phi_\alpha^{-1}(B_{0}(r_{\alpha}R))   \right)\cap \R^{n}_{-}} \left(\Delta^{l/2}\left( \left( B^1_\alpha \circ {\mathcal{T}^{1}_{x_\alpha}}  \right) - \left( B^2_\alpha \circ {\mathcal{T}^{1}_{x_\alpha}}  \right)  \right) \right)^{2}~ dx    \leq \notag\\
& C\sum_{l=0}^{k} \int_{  B_{0}(R)  \cap \R^{n}_{-}}  \left[ \Delta^{l/2}   \left(     \eta^{a}( r _{\alpha}\cdot)u  \right)   -  \Delta^{l/2}  \left( \eta^{b} \left(   \Phi_\alpha ( r _{\alpha}\cdot)\right)    u \left(  r_{\alpha}^{-1}\Phi_\alpha ( r_{\alpha}\cdot) \right) \right) \right]^{2}~ dx =o(1) 
\end{align*}
where $\Phi_\alpha:={\mathcal{T}^{2}_{x_{\alpha}}}^{-1}  \circ{\mathcal{T}^{1}_{x_\alpha}}$ and $d(\Phi_\alpha)_0=Id$. Similarly to the case \eqref{case:1}, we get that
\begin{equation*}
 \lim_{R\to +\infty}\lim_{\alpha\to +\infty}\sum_{l=0}^{k} \int_{M \setminus  D_\alpha(R)  }  \left(\Delta_{g}^{l/2}\left( B^1_\alpha-B^2_\alpha \right) \right)^{2}~ dv_{g} =0.
\end{equation*}
This completes the proof of Proposition \ref{prop:uniqueness}.\qed

\section{Preliminary analysis}\label{sec:prelim}

The proof of Theorem \ref{th:1} goes through four steps. All results are up to a subsequence. We let $(\ua)_\alpha\in H_{k,0}^2(M)$ be a Palais-Smale sequence for $I$. 

\medskip\noindent{\bf Step 1:} We claim that $(\ua)_\alpha$ is bounded in $H_{k,0}^2(M)$.

\smallskip\noindent{\it Proof of the claim:} Since $(\ua)$ is a Palais-Smale sequence, we have that
\begin{eqnarray*}
   \left\langle DI(u_\alpha),u_\alpha\right\rangle &=& 
    \int_M( \Delta_g^{k/2} u_\alpha )^2 ~{dv_g} + 
  \sum_{\alpha =0}^{k-1} \int_{M}A_{l}(\nabla^{l} u_\alpha,\nabla^{l} u_\alpha) ~dv_{g}\\
 && -  \int_M| u_\alpha |^{ 2_k^{\sharp}} ~{dv_g} = o \left(   \left\| u_\alpha \right\|_{H_{k}^2} \right)
 \end{eqnarray*}
Therefore 
\begin{equation}{\label{L^{2*}b'dd}}
     \int_M | u_\alpha |^{ 2_k^{\sharp}}~{dv_g} =\frac{n}{k} I(u_\alpha)
   + o \left(   \left\| u_\alpha \right\|_{H_{k}^2} \right)  \leq C + o \left(   \left\| u_\alpha \right\|_{H_{k}^2} \right)
\end{equation}
Since $(I(\ua))_\alpha$ is bounded, then putting together these equalities yields
\begin{equation*}
  \left\| u_\alpha \right\|_{H_{k}^2}^2 \leq  C+  C \left\| u_\alpha \right\|_{H_{k-1}^2}^2 +C\int_M| u_\alpha |^{ 2_k^{\sharp}} ~{dv_g}
\end{equation*}
Now since the embedding of $H_{k,0}^2(M)$ in $H_{0,k-1}^2(M)$ is compact, then for any $\varepsilon >0 $ there exists a $B_{\varepsilon} >0$ such that $\left\| u \right\|_{H_{k-1}^2}^2 \leq  \varepsilon \left\| u \right\|_{H_{k}^2}^2 +
 B_{\varepsilon} \left\| u \right\|_{{{ 2_k^{\sharp}}} }^2$ for all $u \in H_{k}^2(M)$. Therefore, taking $\varepsilon>0$ small enough, we get that 
\begin{equation*}
  \left\| u_\alpha \right\|_{H_{k}^2}^2 \leq  C+C\int_M| u_\alpha |^{ 2_k^{\sharp}} ~{dv_g}
\end{equation*}
Then using \eqref{L^{2*}b'dd} we get that $\left\| u_\alpha \right\|_{H_{k}^2}^2 \leq C+C\left\| u_\alpha \right\|_{H_{k}^2} $ for all $\alpha$, and therefore the sequence $(u_\alpha)$ is bounded in  $H_{k,0}^2(M)$. This proves the claim.\qed

\medskip\noindent  Since $(u_\alpha)$ is bounded in  $H_{k,0}^2(M)$, there exists  $u_\infty \in H_{k,0}^2(M)$ such that 
 \begin{eqnarray}{\label{weakly and strongly}}
          \left \{ \begin{array} {lc}
                  u_\alpha \rightharpoonup u_\infty \quad \text{weakly in } H_{k,0}^2(M)\hbox{ and } L^{2_k^{\sharp}}(M),\\
                  u_\alpha \rightarrow u_\infty \quad \text{strongly in } H_{l,0}^2(M) \hbox{ and in }L^{q}(M) \hbox{ for } l<k,\,q< 2_k^{\sharp}, \\
                                   u_\alpha(x) \rightarrow u_\infty(x) \quad \text{a.e in } M
           \end{array} \right.
\end{eqnarray}
We define $v_\alpha := u_\alpha- u_\infty $.

\medskip\noindent{\bf Step 2:} We claim that
\begin{enumerate}
 \item  $DI(u_\infty) =0$
 \item  $(v_\alpha)$ is  a Palais-Smale sequence for the functional $J$ on the space $H_{k,0}^2(M)$,
 \item $ J(v_\alpha) = I(u_\alpha) - I(u_\infty) + o(1)$ as $\alpha\to +\infty$.
\end{enumerate}
where
\begin{equation*}
  J(u):= \frac{1}{2}\int_M( \Delta_g^{k/2} u)^2~{dv_g} 
 - \frac{1}{2_k^{\sharp}} \int_M| u |^{ 2_k^{\sharp}}~{dv_g}\hbox{ for }u\in H_{k,0}^2(M)
 \end{equation*}
 \smallskip\noindent{\it Proof of the claim:} We fix $\varphi \in H_{k,0}^2(M)$. We have that
\begin{align}\label{eq:crit:u}
 \left\langle DI(u_\alpha),\varphi\right\rangle =  &
  \int_M\Delta_g^{k/2} u_\alpha \Delta_g^{k/2} \varphi  ~{dv_g} + 
  \sum_{\alpha =0}^{k-1} \int_M
   A_{l}(g) (\nabla^{l} u_\alpha,\nabla^{l} \varphi) \notag\\
     & - \int_M \left| u_\alpha \right|^{ 2_k^{\sharp}-2}u_\alpha  \varphi ~{dv_g} = o(1)
\end{align}
The following classical integration Lemma will be often used in the sequel (see Lemma 6.2.7 in Hebey \cite{hebey.f.book} for a proof):
\begin{lemma}\label{int.lemma} Let $(M,g)$ be a Riemannian manifold. If $(f_{\alpha})$ is a bounded sequence in $L^{p}(M)$, $1< p< + \infty$,  such that $f_{\alpha} \rightarrow f $ $a.e$ in $M$, then $f \in L^{p}(M)$ and $f_{\alpha} \rightharpoonup f$ weakly in $L^{p}(M)$. 
\end{lemma}
\noindent Since $(\left| u_\alpha \right|^{ 2_k^{\sharp}-2}u_\alpha)_\alpha$ is  bounded in $L^{\frac{\crit}{\crit-1}}$ and converges a.e., Lemma \ref{int.lemma} yields
\begin{equation}\label{lim:u0}
 \int_M \left| u_\alpha \right|^{ 2_k^{\sharp}-2}u_\alpha  \varphi ~{dv_g} =
 \int_M \left| u_\infty \right|^{ 2_k^{\sharp}-2}u_\infty  \varphi ~{dv_g} + o(1)
\end{equation}
Therefore, the weak convergence of $(u_\alpha)$ to $u_\infty$, \eqref{eq:crit:u} and \eqref{lim:u0} yield that $u_\infty$ is a weak solution to \eqref{eqn1}. This proves point (1) of Step 2.

\medskip\noindent We now estimate $I(\ua)$. From \eqref{weakly and strongly} we have
\begin{equation*}
    \int_M( \Delta_g^{k/2} u_\alpha )^2 ~{dv_g} - \int_M( \Delta_g^{k/2} u_\infty )^2 ~{dv_g}
    = \int_M( \Delta_g^{k/2} v_\alpha )^2 ~{dv_g}   + o(1) ,
  \end{equation*}
   \begin{equation*}
  \sum_{l =0}^{k-1}  \int_MA_{l} (\nabla^{l} u_\alpha,\nabla^{l} u_\alpha) ~dv_{g} =
   \sum_{l =0}^{k-1}  \int_MA_{l} (\nabla^{l} u_\infty,\nabla^{l} u_\infty) ~dv_{g} + o(1)
 \end{equation*}
 The following two inequalities will be of constant use in the sequel: for any $1<p<+\infty$, there exists $C>0$ such that
\begin{equation}\label{ineq:useful:2}
\left| \ |a+b|^p - |a|^p  -|b|^p   \  \right| \leq C \left( |a|^{p-1}|b|+ |b|^{p-1} |a|\right) 
\end{equation}
\begin{equation}\label{ineq:useful}
\left| \ |a+b|^p(a+b) - |a|^p a -|b|^p b  \  \right| \leq C \left( |a|^p|b|+ |b|^p |a|\right) 
\end{equation}
for all $a,b\in\R$. It then follows from \eqref{ineq:useful:2} that
$$\left| |u_\alpha|^{\crit}-|u_\infty|^{\crit}-|v_\alpha|^{\crit}\right |\leq C\left(|v_\alpha|^{\crit-1}|u_\infty|+|u_\infty|^{\crit-1}|v_\alpha|\right),$$
and then using Lemma \ref{int.lemma}, we get that
 \begin{equation*}
  \int_M| u_\alpha |^{ 2_k^{\sharp}} ~{dv_g} -  \int_M| u_\infty |^{ 2_k^{\sharp}} ~{dv_g}
  = \int_M| v_\alpha |^{ 2_k^{\sharp}} ~{dv_g} + o(1)
\end{equation*}
Hence $ I(u_\alpha) - I(u_\infty) = J(v_\alpha)  + o(1)$ as $\alpha\to +\infty$, which proves point (3) of Step 2. 

\medskip\noindent Next we show the sequence $(v_\alpha)$ is  a Palais-Smale sequence for the functional $J$ on $H_{k,0}^2(M)$. Let
$\varphi \in H_{k,0}^2(M)$, we have
\begin{equation}\label{crit:V:U}
 \left\langle DJ(v_\alpha),\varphi\right\rangle =
 \left\langle DI(u_\alpha),\varphi\right\rangle - \left\langle DI(u_\infty),\varphi\right\rangle + 
 \int_M \varPhi_\alpha \varphi ~{dv_g}   + o(\left\|\varphi \right\|_{H_k^2})
\end{equation}
where
\begin{equation*}
 \varPhi_\alpha : =  \left| v_\alpha + u_\infty \right|^{ 2_k^{\sharp}-2}(v_\alpha +u_\infty)    -  \left| u_\infty \right|^{ 2_k^{\sharp}-2}u_\infty   -
  \left| v_\alpha \right|^{ 2_k^{\sharp}-2}v_\alpha  
\end{equation*}
Inequality \eqref{ineq:useful} and  H\"{o}lder's inequality yield
\begin{equation}\label{bnd:phi}
 \left|  \int_M \varPhi_\alpha \varphi ~{dv_g} \right|  \leq
 C \left( \left\| |v_\alpha|^{{2_k^{\sharp} -2 }} u_\infty \right\|_{\frac{2_k^{\sharp}} {2_k^{\sharp} -1 }} + 
 \left\||u_\infty|^{{2_k^{\sharp} -2 }} v_\alpha \right\|_{\frac{2_k^{\sharp}} {2_k^{\sharp} -1 }} \right)
 \left\| \varphi \right\|_{2_k^{\sharp}} 
\end{equation}
Since $v_\alpha \rightharpoonup 0 \quad \text{in } L^{2_k^{\sharp}}(M)$, Lemma \ref{int.lemma} yields
\begin{equation*}
 \left\| |v_\alpha|^{{2_k^{\sharp} -2 }} u_\infty \right\|_{\frac{2_k^{\sharp}} {2_k^{\sharp} -1 }} + 
 \left\||u_\infty|^{{2_k^{\sharp} -2 }} v_\alpha \right\|_{\frac{2_k^{\sharp}} {2_k^{\sharp} -1 }} = o(1)
\end{equation*}
Since $(u_\alpha)$ is a Palais-Smale for $I$, then \eqref{crit:V:U}, \eqref{bnd:phi} and the continuous embedding $H_{k,0}^2 (M)\hookrightarrow L^{2_k^{\sharp}}(M)$ yields $ \left\langle DJ(v_\alpha),\varphi\right\rangle = o(\left\|\varphi \right\|_{H_k^2})$ as $\alpha\to +\infty$ uniformly wrt $\varphi\in H_{k,0}^2(M)$. This proves the claim and ends Step 2.\qed

\medskip\noindent The next lemma adresses the compactness of a Palais-Smale sequence for small energy. It will be generalized to the case of small local energy in Proposition \ref{prop:fund}.

\medskip\noindent {\bf Step 3:} Let $(v_\alpha)$ be a Palais-Smale sequence for $J$ on $H_{k,0}^2(M)$. We assume that $ v_\alpha \rightharpoonup  0 $  weakly in
 $H_{k,0}^2(M)$, and  that $J(v_\alpha) \rightarrow \beta$ with $\beta <  \beta^\sharp
$, where $\beta^\sharp$ is as in \eqref{def:beta}. We claim that $v_\alpha \rightarrow 0 $ strongly
 in $H_{k,0}^2(M)$.

\medskip\noindent{\it Proof of the claim:} Since $(v_\alpha)$ is bounded and $\langle DJ(v_\alpha),v_\alpha\rangle=o(\Vert v_\alpha\Vert_{H_k^2})$, we get that
 \begin{equation}{\label{2.6}}
  J(v_\alpha)  = \frac{k}{n} \int_M( \Delta_g^{k/2} v_\alpha)^2 ~{dv_g} + o(1)
  = \frac{k}{n}  \int_M| v_\alpha |^{ 2_k^{\sharp}} ~{dv_g} + o(1) = \beta + o(1).
 \end{equation}
As a consequence, $\beta \geq 0$. It follows from Mazumdar \cite{mazumdar.gjms} that for any $\varepsilon >0 $ there exists  $B_{\varepsilon} >0$ such that 
 \begin{equation}\label{ineq:sobo}
  \left\| u\right\|_{\crit}^2 \leq \left(K_0(n,k) + \varepsilon \right) \int_{\tilde{M}} ( \Delta_{g}^{k/2} u )^2 ~{dv_{g}} + 
 B_{\varepsilon}\left\| u\right\|_{H_{k-1}^2}^2 
 \end{equation}
for all $u\in H_{k}^2(\tilde{M})$. Applying this inequality to $v_\alpha$, the strong convergence to $0$ in $H_{k-1}^2$ and \eqref{2.6} yield 

\begin{equation*}
 \left(  \frac{n}{k} \beta \right)^{2/{2_k^{\sharp}}} \leq  \left(K_0(n,k) + \varepsilon \right) \frac{n}{k} \beta
\end{equation*}
Letting $\varepsilon \to 0$ and using $0\leq \beta < \beta^\sharp$, we get that $\beta=0$, and then \eqref{2.6} yields $v_\alpha\to 0$ strongly in $H_{k,0}^2(M)$. This proves the claim and ends Step 3.\qed

\medskip\noindent{\bf Step 4: Proof of Theorem \ref{th:1}.} Let $(u_\alpha)$ be a Palais-Smale sequence for the functional $I$ on the space $H_{k,0}^2(M)$. By substracting the weak limit $u_\infty$, we get a Palais-Smale sequence $(v_\alpha)$ for the functional $J$ with energy $ J(v_\alpha) = I(u_\alpha)- I(u_\infty) + o(1)$ as $\alpha\to +\infty$. If $v_\alpha\to 0$ strongly in $H_{k,0}^2(M)$, then we end the process. If not, we apply Lemma \ref{bubble.lemma} to substract a bubble modeled on $v\in{\mathcal D}_k^2(\rn)\setminus \{0\}$ and we get a new Palais-Smale sequence for $J$, but with the energy decreased by $E(v)$. If the resulting sequence goes strongly to $0$, we stop the process, if not, we iterate it again. This process must stop since the energy $E(v)\geq\beta^\sharp$ and after finitely many steps, the energy goes below the critical threshold $\beta^\sharp$ and then the convergence is strong by Step 3. This proves Theorem \ref{th:1}.\qed

\medskip\noindent The rest of the paper is devoted to the proof of Lemma \ref{bubble.lemma}.

 \section{Extraction of a Bubble}\label{sec:extrac}
 In the sequel, for any $(M,g)$ as in the introduction, we let $H_k^2(M)$ be the completion of $\{u\in C^\infty(M): \Vert u\Vert_{H_k^2}<+\infty\}$ for the norm $\Vert\cdot\Vert_{H_k^2}$. The space $H_{k,0}^2(M)$ is then a closed subspace of $H_k^2(M)$. The following lemma is the main ingredient in the proof of Theorem \ref{th:1} 
\begin{lemma}{\label{bubble.lemma}}
 Let $(v_\alpha)$ be  a  Palais-Smale sequence for the functional $J$ on  $H_{k,0}^2(M)$  such that $v_\alpha \rightharpoonup 0$  weakly  in $H_{k,0}^2(M)$ but not strongly. Then there exists a bubble $(B_{x_\alpha,r_\alpha}(v))$   such that upto a subsequence, the following holds:
 \begin{itemize}
 \item $w_\alpha:=v_\alpha-B_{x_\alpha,r_\alpha}(v)$ is a Palais-Smale sequence for $J$,
 \item $J(w_\alpha)=J(v_\alpha)-E(v)+o(1)$ as $\alpha\to +\infty$.
 \end{itemize}
 \end{lemma}

\noindent The proof of this lemma goes through 10 steps.

\medskip\noindent{\bf Step 1:} We prove a strong convergence Lemma for small energies. This is a localized version of Step 3 of Section \ref{sec:prelim}.

 \begin{proposition}\label{prop:fund}
Let $(N,g_\infty)$ be a Riemannian manifold with positive injectivity radius. 

\smallskip\noindent$\bullet$ Let $(g_i)_i$ be metrics on $N$ such that $g_i\to g_\infty$ in $C^p_{loc}$ as $i\to +\infty$ for all $p$.

\smallskip\noindent$\bullet$ Let $(P_i)_i$ be a family of operators on $C^\infty(N)$ such that 
$$P_i:=\Delta_{g_i}^k+\sum_{l=0}^{k-1}(-1)^l\nabla^{i_1...i_l}\left((A^{i}_l)_{i_1...i_lj_1...j_l}\nabla^{j_1...j_l}\right)$$ 
with families of symmetric tensors $(A_l^i)\to A_l$ in $C^p_{loc}$ as $i\to +\infty$ for all $p$.

\smallskip\noindent$\bullet$ We fix $\Omega\subset N$ an open smooth domain, and we define
\begin{equation}\label{hyp:j:prime}
J_i(u):=\frac{1}{2}\int_\Omega u P_i u\, dv_{g_i}-\frac{1}{\crit}\int_\Omega |u|^{\crit}\, dv_{g_i}\hbox{ for }u\in H_{k}^2(\Omega),
\end{equation}
such that $J_i$ is $C^1$. Here, the background metric is $g_\infty$.

\smallskip\noindent$\bullet$ We let $(u_i)\in H_{k,0}^2(\Omega)$ and $u_\infty\in H_{k,0}^2(\Omega)$ be such that $u_i\rightharpoonup u_\infty$ weakly in $H_{k,0}^2(\Omega)$ as $i\to +\infty$.

\smallskip\noindent$\bullet$ We assume that there exist a compact $K\subset N$ such that
$$\lim_{i\to +\infty}\sup_{u\in H_{k,0}^2(\Omega),\, \hbox{Supp }\varphi\subset K}\frac{\langle DJ_i(u_i),\varphi\rangle}{\Vert \varphi\Vert_{H_k^2(\Omega)}}=0$$

\smallskip\noindent$\bullet$ We assume that there exists $K_\infty>0$ and $C\geq 0$ such that 
\begin{equation}\label{ineq:sobo:2}
\left(\int_N |u|^{\crit}\, dv_{g_\infty}\right)^{\frac{2}{\crit}}\leq K_\infty \int_N (\Delta_{g_\infty}^{k/2} u)^2\, dv_{g_\infty}+ C\Vert u\Vert_{H_{k-1}^2}^2\hbox{ for all }u\in C^\infty_c(N).
\end{equation}
We fix $x_0\in \Omega$ and $\delta\in (0, i_{g_\infty}(N)/2)$. We assume that
\begin{equation}\label{hyp:lemma}\left\{\begin{array}{l}
B_{x_0} (2\delta) \subset K \;\hbox{$($the ball is wrt $g_\infty)$},\\
\displaystyle{\int_{B_{x_0} (2\delta)  \cap\Omega}|u_i|^{\crit}\, dv_{g_i}\leq \left( \frac{1}{2 K_\infty } \right)^{\frac{\crit}{\crit-2}}}\hbox{ for all }i\in\N.
\end{array}\right.
\end{equation}
Then $u_i\to u_\infty$ strongly in $H_k^2(B_{ x_0 }(\delta)\cap\Omega)$.
\end{proposition}

\medskip\noindent{\it Proof of Proposition \ref{prop:fund}:} Up to extracting a subsequence, we assume that $u_i\to u_\infty$ strongly in $H_{k-1}^2(\omega)$ as $i\to +\infty$ for $\omega\subset \Omega$ relatively compact and $u_i(x)\to u_\infty(x)$ as $i\to +\infty$ for a.e. $x\in \Omega$. 
Let $\eta\in C^\infty(N)$ such that $\eta(x)=1$ for $x\in B_{x_0}(\delta)$ and $\eta(x)=0$ for $x\in N\setminus B_{x_{0}}(2\delta)$. Since $\eta$ has compact support, we get that $\eta^2(u_i-u_\infty)\in H_{k,0}^2(\Omega)$ is uniformly bounded in $H_{k,0}^2(\Omega)$. Since $B_{ x_0 }(2\delta)\subset K$, it then follows from hypothesis \eqref{hyp:j:prime} that
$$\langle DJ_i(u_i),\eta^2(u_i-u_\infty)\rangle=o(1)\hbox{ as }i\to +\infty.$$
Since $\eta^2(u_i-u_\infty)\to 0$ strongly in $H_{k-1}^2(\Omega)$, we then get that
\begin{equation}\label{eq:1}
\int_\Omega \Delta_{g_i}^{k/2}u_i\Delta_{g_i}^{k/2}(\eta^2(u_i-u_\infty))\, dv_{g_i}=\int_\Omega |u_i|^{\crit-2}u_i\eta^2(u_i-u_\infty)\, dv_{g_i}+o(1)
\end{equation}
as $i\to +\infty$. The weak convergence of $u_i$ to $u_\infty$ and the strong convergence of $g_i$ to $g_\infty$ on compact sets yields 
\begin{equation}\label{eq:2}
\int_\Omega \Delta_{g_i}^{k/2}u_i\Delta_{g_i}^{k/2}(\eta^2(u_i-u_\infty))\, dv_{g_i}=\int_\Omega \Delta_{g_i}^{k/2}(u_i-u_\infty)\Delta_{g_i}^{k/2}(\eta^2(u_i-u_\infty))\, dv_{g_i}+o(1)
\end{equation}
as $i\to +\infty$. As one checks, for any $\varphi\in H_k^2(\Omega)$, we have that $\Delta_{g_i}^{k/2}\varphi\Delta_{g_i}^{k/2}(\eta^2\varphi)=\left(\Delta_{g_i}^{k/2}(\eta\varphi)\right)^2+\sum_{p<k,l\leq k}\nabla^p \varphi\star\nabla^l\varphi,$ where $A\star B$ denotes a linear combination of bilinear forms in $A$ and $B$. Therefore, using again the strong convergence of $\eta^2(u_i-u_\infty)$ to $0$ in $H_{k-1}^2$, we get that 
\begin{equation}\label{eq:3}
\int_\Omega \Delta_{g_i}^{k/2}u_i\Delta_{g_i}^{k/2}(\eta^2(u_i-u_\infty))\, dv_{g_i}=\int_\Omega \left(\Delta_{g_i}^{k/2}(\eta(u_i-u_\infty))\right)^2\, dv_{g_i}+o(1)
\end{equation}
as $i\to +\infty$. Moreover, since $|u_i|^{\crit-2}\eta^2(u_i-u_\infty)$ is uniformly bounded in $L^{\crit/(\crit-1)}$ and goes to $0$ almost everywhere as $i\to +\infty$, then it goes weakly to $0$ in $L^{\crit/(\crit-1)}$, and then $\int_\Omega |u_i|^{\crit-2}\eta^2(u_i-u_\infty)u_\infty\, dv_{g_i}\to 0$ as $i\to +\infty$. Therefore, plugging \eqref{eq:2} and \eqref{eq:3} into \eqref{eq:1}, we get that 
$$\int_\Omega \left(\Delta_{g_i}^{k/2}(\eta(u_i-u_\infty))\right)^2\, dv_{g_i}=\int_\Omega |u_i|^{\crit-2}(\eta(u_i-u_\infty))^2\, dv_{g_i}+o(1)$$
as $i\to +\infty$. Since $g_i\to g_\infty$ as $i\to +\infty$ in $C^p$ locally on compact sets and $\eta(u_i-u_\infty)$ is uniformly bounded in $H_k^2(\Omega)$, we get that
$$\int_\Omega \left(\Delta_{g_\infty}^{k/2}(\eta(u_i-u_\infty))\right)^2\, dv_{g_\infty}=\int_\Omega |u_i|^{\crit-2}(\eta(u_i-u_\infty))^2\, dv_{g_\infty}+o(1)$$
as $i\to +\infty$. H\"older's inequality, the Sobolev inequality \eqref{ineq:sobo:2}, the convergence of $(g_i)$, the strong convergence in $H_{k-1}^2$ and \eqref{hyp:lemma} then yields
\begin{eqnarray*}
&&\int_\Omega \left(\Delta_{g_\infty}^{k/2}(\eta(u_i-u_\infty))\right)^2\, dv_{g_\infty}\\
&&\leq \left(\int_{B_{x_0}(2\delta) \cap\Omega} |u_i|^{\crit}\, dv_{g_\infty}\right)^{\frac{\crit-2}{\crit}}\left(\int_N|\eta(u_i-u_\infty)|^{\crit}\, dv_{g_\infty}\right)^{\frac{2}{\crit}}+o(1)\\
&&\leq\frac{1}{2K_\infty}\left(K_\infty \int_N \left(\Delta_{g_\infty}^{k/2}(\eta(u_i-u_\infty))\right)^2\, dv_{g_\infty}+C\Vert \eta(u_i-u_\infty)\Vert^2_{H_{k-1}^2}\right)+o(1)\\
&&\leq  \left(\int_{B_{x_0}(2\delta)\cap\Omega} |u_i|^{\crit}\, dv_{g_i}\right)^{\frac{\crit-2}{\crit}}K_\infty \int_\Omega \left(\Delta_{g_\infty}^{k/2}(\eta(u_i-u_\infty))\right)^2\, dv_{g_\infty}+o(1)\\
\end{eqnarray*}
as $i\to +\infty$. Therefore, we get that $\Vert\Delta_{g_\infty}^{k/2}(\eta(u_i-u_\infty))\Vert_2\to 0$ as $i\to +\infty$. Since $\eta(u_i-u_\infty)\to 0$ strongly in $H_{k-1}^2$ and $\eta$ has compact support, we get that  $\eta(u_i-u_\infty)\to 0$ strongly in $H_k^2(\Omega)$, and therefore $u_i\to u_\infty$ in $H_k^2(B_{x_0}(\delta)\cap\Omega)$. Note that this is up to a subsequence. Indeed, by uniqueness, the convergence holds for the initial sequence $(u_i)$. This proves Proposition \ref{prop:fund}. \hfill$\Box$

\medskip\noindent{\bf Step 2:}  Since 
   $\left\langle DJ(v_\alpha),v_\alpha\right\rangle = o(1) $, one has
  \begin{equation*}
  J(v_\alpha)= \frac{k}{n}\int_{M} |v_\alpha|^{\crit} ~{dv_g} + o(1) =\beta+o(1)~ \text{as}\ \ \alpha \rightarrow + \infty
 \end{equation*}
where $\beta:=\lim_{\alpha\to +\infty}J(v_\alpha)$. By Step 3 of Section \ref{sec:prelim},  $\beta \geq \beta^{\sharp}$. Therefore, since $\overline{M}$ is compact, for any $r_0 >0$, there exists $y_{0} \in \overline{M}$ and $\lambda_0 >0 $ such that 
\begin{equation*}
 \int_{B_{y_{0}}(r_0)\cap M } |v_\alpha|^{\crit} ~ {dv_g} \geq \lambda_0
\end{equation*}
\medskip\noindent For any $r>0$, we set 
\begin{equation}\label{def:mu}
 \mu_\alpha (r):= \max_{x \in \overline{M}} \int_{B_{x}(r)\cap M} |v_\alpha|^{\crit}  ~{dv_g},
\end{equation}
the Levy concentration function. In particular,  $\mu_\alpha (r_0) \geq \lambda_0 $ for all $\alpha$. We fix 
$$0<\lambda<\epsilon_0:=\min\left\{\lambda_0,\frac{1}{(2 K_0(n,k))^{\crit/(\crit-2)}}\right\}$$
where $K_0(n,k)$ is the best constant in the Euclidean Sobolev inequality \eqref{sinq}. Since $\mu_\alpha(0)=0$, there exists  $(r_\alpha)_\alpha \in (0, r_0)$ and $(x_\alpha)_\alpha \in \overline{M}$ such that:
\begin{equation}\label{def:ra:xa}
 \lambda = \mu_\alpha (r_\alpha) = \int_{B_{x_\alpha}(r_\alpha)\cap M} |v_\alpha|^{\crit}~ {dv_g}
\end{equation}
\medskip

\medskip\noindent {\bf Step 3:} We claim that $\lim_{\alpha\to +\infty}r_\alpha=0$.

\smallskip\noindent{\it Proof of the claim.} We argue by contradiction. If $(r_\alpha)$ does not go to $0$ up to a subsequence, we get that there exists $\delta\in (0, i_g(\tilde{M})/2)$ such that for all $x\in M$, we have that $\int_{B_{x}(2\delta)\cap M}|v_\alpha|^{\crit}\, dv_g\leq \lambda$ for all $\alpha$. We apply Proposition \ref{prop:fund} with $(N,g_\infty)=(\tilde{M}, g)$, $\Omega=M$, $P_\alpha=P$, $g_\alpha=g$, $J_\alpha=J$, and the Sobolev inequality \eqref{ineq:sobo} of \cite{mazumdar.gjms}, and we get $v_\alpha\to 0$ as $\alpha\to +\infty$ in $H_k^2(M\cap B_{x}(\delta))$ for all $x\in M$. With a finite covering, we get that $v_\alpha\to 0$ as $\alpha\to +\infty$ strongly in $H_{k,0}^2(M)$, contradicting our initial hypothesis. This proves the claim and ends Step 3.\qed

\medskip\noindent  First  assume that
\begin{equation}\label{hyp:A}
\lim_{\alpha\to +\infty}\frac{d(x_\alpha, \partial M)}{r_\alpha}=+\infty.
\end{equation} 
We define
$$\tva(x):=r_\alpha^{\frac{n-2k}{2}}\ua (exp_{\xm}(r_\alpha x))\hbox{ for }|x|<\frac{i_g(\tilde{M})}{r_\alpha}\hbox{ and }|x|<\frac{d(\xm,\partial M)}{r_\alpha}$$

\medskip\noindent {\bf Step 4:} Suppose  that \eqref{hyp:A} holds. We claim that there exists $v\in {\mathcal D}_k^2(\rn)$ such that for any $\eta\in C^\infty_c(\rn)$, we have that
\begin{equation*}
\eta \tva\rightharpoonup \eta v\hbox{ weakly in }{\mathcal D}_k^2(\rn)\hbox{ as }k\to +\infty.
\end{equation*}

\smallskip\noindent{\it Proof of the claim.} Fix $\eta\in C^\infty_c(\rn)$, and let $R_0>0$ be such that $\hbox{Supp }\eta\subset B_{0}(R_{0})$. We define
$$\eta_\alpha(x):=\eta\left(\frac{\hbox{exp}_{\xm}^{-1}(x)}{r_\alpha}\right)\hbox{ for }x\in B_{\xm }(R_0r_\alpha),\hbox{ and }\eta_\alpha(x):=0\hbox{ outside.}$$
Up to a subsequence, there exists $x_0\in \tilde{M}$ and $\tau>0$ such that $B_{\xm }(R_0r_\alpha) \subset B_{x_0}(\tau)\subset \tilde{M}$. It then follows from the comparison Lemma  9.1  of Mazumdar \cite{mazumdar.gjms} that there exists $C>0$ such that 
\begin{equation*}
\int_{B_{0}(R_0r_\alpha)}\left(\Delta^{k/2}[(\eta_\alpha v_\alpha)\circ \hbox{exp}_{\xm}]\right)^2\, dx\leq C \int_{B_{\xm }(R_0r_\alpha)} \left(\Delta_g^{k/2}(\eta_\alpha v_\alpha)\right)^2\, dv_g 
\end{equation*}
for all $\alpha$. With a change of variable, rough estimates of the differential terms and H\"older's inequality, we then get 
\begin{eqnarray}
&&\int_{B_{0}(R_{0})}\left(\Delta^{k/2}(\eta \tva)\right)^2\, dx\leq C \sum_{l=0}^k\int_{B_{\xm }(R_0r_\alpha)} |\nabla^l\um|_g^2|\nabla^{k-l}\eta_\alpha|_g^2\, dv_g\nonumber\\
&&\leq   C\sum_{l=0}^{k}\int_{B_{\xm }(R_0r_\alpha)} r_\alpha^{2(l-k)}|\nabla^lv_\alpha|_g^2\, dv_g\leq   C\sum_{l=0}^{k}\Vert\nabla^lv_\alpha\Vert_{\frac{2n}{n-2(k-l)}}^2\label{eq:10}
\end{eqnarray}
It follows from Sobolev's embedding theorem that $H_{k-l}^{2}(M)\subset L^{\frac{2n}{n-2(k-l)}}(M)$ for all $l=0,...,k$ and that this embedding is continuous. Since $(v_\alpha)_\alpha$ is bounded in $H_k^2$, then $(\nabla^lv_\alpha)_\alpha$ is uniformly bounded in $H_{k-l}^2$ (with tensorial values), and then there exists $C>0$ such that 
\begin{equation}\label{ineq:nabla}
\Vert\nabla^lv_\alpha\Vert_{\frac{2n}{n-2(k-l)}}\leq C\Vert v_\alpha\Vert_{H_k^2}\leq C'
\end{equation}
for all $\alpha>0$ and $l=0,...,k$. It then follows from \eqref{eq:10} that $(\eta\tva)_\alpha$ is bounded in ${\mathcal D}_k^2(\rn)$. Therefore, up to a subsequence, there exists $v_\eta\in {\mathcal D}_k^2(\rn)$ such that $\eta\tva\rightharpoonup v_\eta$ weakly in ${\mathcal D}_k^2(\rn)$ as $\alpha\to +\infty$. A classical diagonal argument then yields the existence $v\in H_{k,loc}^2(\rn)$ such that $\eta\tva\rightharpoonup \eta v$ weakly in ${\mathcal D}_k^2(\rn)$ as $\alpha\to +\infty$.  We fix $R>0$. For any $R'>R$, a change of variables and \eqref{ineq:nabla} yields
$$\int_{B_{0}(R)} |\nabla^l\eta_{R'}\tva|_{g_\alpha}^{\frac{2n}{n-2(k-l)}}\, dv_{g_\alpha}\leq \int_{B_{\xm }(R_0r_\alpha)} |\nabla^lv_\alpha|_g^{\frac{2n}{n-2(k-l)}}\, dv_g\leq C$$
where $g_\alpha:=\hbox{exp}_{\xm}^\star g(r_\alpha \cdot)$. Using weak convergence and convexity, letting $\alpha\to +\infty$ and then $R\to +\infty$ yields $|\nabla^lv|\in L^{\frac{2n}{n-2(k-l)}}(\rn)$. As one checks, we then have that the sequence  $(\eta_Rv)_R$ is a Cauchy sequence in ${\mathcal D}_k^2(\rn)$, and then we get that $v\in {\mathcal D}_k^2(\rn)$. This ends the proof of the claim,  and ends Step 4.\qed

\medskip
\noindent {\bf Step 5:} {\label{strong conv}}
 We assume that \eqref{hyp:A} holds. We let $v\in {\mathcal D}_k^2(\rn)$ as in Claim 3. We claim that $v\not\equiv 0$ is a weak solution to $\Delta^kv=|v|^{\crit-2}v$ in ${\mathcal D}_k^2(\rn)$. 

\medskip\noindent{\it Proof of the claim}. We fix $R>0$ and we apply Proposition \ref{prop:fund} with $(N, g_\infty):=(\rn, \hbox{Eucl})$ and $\Omega:=\rn$. As above, we define a family of smooth metrics $(g_\alpha)_\alpha$ such that $g_\alpha(x):=\hbox{exp}_{\xm}^\star g(r_\alpha x)$ for $x\in B_0(3R)$, $g_\alpha(x)=\hbox{Eucl}$ for $x\in\rn\setminus B_{0}(4R)$, and $g_\alpha\to \hbox{Eucl}$ in $C^p_{loc}(\rn)$ as $\alpha\to +\infty$ for all $p$. Let $\varphi\in C^\infty_c(\rn)$ be such that $\hbox{Supp }\varphi\subset B_0(R)$. We define 
$$\varphi_\alpha(x):=r_\alpha^{-\frac{n-2k}{2}}\varphi\left(\frac{\hbox{exp}_{\xm}^{-1}(x)}{r_\alpha}\right)$$
for all $x\in M$. As one checks, $\varphi_\alpha$ is well-defined and has support in $B_{\xm }(Rr_\alpha)$. Moreover, using the comparison Lemma 9.1  in Mazumdar \cite{mazumdar.gjms} and arguing as in Step 4, we get that $\Vert\varphi_\alpha\Vert_{H_{k,0}^2(M)}\leq C(R)\Vert\varphi\Vert_{H_{k,0}^2(\rn)}$ for all $\alpha>0$. Since $(\um)$ is a Palais-Smale sequence, we have that
$$\langle DJ(v_\alpha),\varphi_\alpha\rangle=o(\Vert\varphi_\alpha\Vert_{H_{k,0}^2})=o(\Vert\varphi\Vert_{H_{k,0}^2(\rn)})$$
as $\alpha\to +\infty$ uniformly for all $\varphi\in C^\infty_c(\rn)$ such that $\hbox{Supp }\varphi\subset B_0(R)$. With a change of variable, we get  $\langle DJ(v_\alpha),\varphi_\alpha\rangle=\langle DJ_\alpha(\eta_R\tva),\varphi\rangle$ where 
$$J_\alpha(u):=\frac{1}{2}\int_{\rn}(\Delta_{g_\alpha}^{k/2} u)^2\, dv_{g_\alpha}-\frac{1}{\crit}\int_{\rn}|u|^{\crit}\, dv_{g_\alpha}$$
for all $u\in H_k^2(\rn)$. Therefore, $\langle DJ_\alpha(\eta_R\tva),\varphi\rangle=o(\Vert\varphi\Vert_{H_{k,0}^2(\rn)})$ as $\alpha\to +\infty$ uniformly for all $\varphi\in C^\infty_c(\rn)$ such that $\hbox{Supp }\varphi\subset B_R(0)$.

\smallskip\noindent We fix $x_0\in\rn$ such that $B_{x_{0}}(1/2)\subset B_0(R)$. A change of variable yields
$$\int_{B_{x_{0}}(1/2)\cap B_{0}(2R)}|\eta_R\tva|^{\crit}\, dv_{g_\alpha}=\int_{\hbox{exp}_{\xm}(r_\alpha B_{x_{0}}(1/2))}|\um|^{\crit}\, dv_g.$$
For $\alpha>0$ large enough, we have that $\hbox{exp}_{\xm}(r_\alpha B_{x_{0}}(1/2))\subset B_{\hbox{exp}_{\xm}(x_0)}(r_\alpha)$. Therefore, it follows from the definition of $\mu_\alpha$ that
$$\int_{B_{x_{0}}(1/2)\cap B_{0}(2R)}|\eta_R\tva|^{\crit}\, dv_{g_\alpha}\leq \mu_\alpha(r_\alpha)=\lambda<\epsilon_0$$
for all $\alpha$ large enough and $x_0\in\rn$ such that $1/2+|x_0|<R$. With the Sobolev inequality \eqref{sinq} on $\rn$, we apply Proposition \ref{prop:fund} to $(\eta_R\tva)_\alpha$, and we get that 
$$\lim_{\alpha\to +\infty}\eta_R\tva= \eta_Rv\hbox{ strongly in }H_k^2(B_{x_{0}}(1/4)).$$
Using a finite covering, we then have $\tilde{v}_\alpha\to v$ strongly in $H_k^2(B_{0}(R/2))$ as $\alpha\to +\infty$. Sobolev's embedding theorem yield the convergence in $L^{\crit}(B_0(1))$. Since
$$\int_{B_0(1)}|\tva|^{\crit}\, dv_{g_\alpha}=\int_{B_{\xm }(r_\alpha)}|v_\alpha|^{\crit}\, dv_g=\mu_\alpha(r_\alpha)=\lambda>0,$$
passing to the limit $\alpha\to +\infty$ yields $\int_{B_0(1)}|v|^{\crit}\, dx=\lambda\neq 0$, and therefore $v\not\equiv 0$. This proves the claim and ends Step 5.\qed

\medskip\noindent Note that indeed, we have proved that
\begin{equation}\label{strong:cv}
\lim_{\alpha\to +\infty}\tva=v\hbox{ strongly in }H_k^2(B_0(R))\hbox{ for all }R>0.
\end{equation}
We choose a sequence $(\tilde{r}_{\alpha})$ of positive real numbers as in \eqref{def:rt} with $\eta\in C^\infty_c(B_0(\delta))$ (with $\delta\in (0, i_g(\tilde{M}))$) identically $1$ around $0$. As in Definition \ref{def:bubbles}, we set

\begin{align*}
 V_\alpha(x):= B_{x_\alpha,r_\alpha}(v):=\eta\left(\frac{exp_{x_\alpha}^{-1}(x)}{\tilde{r}_\alpha}\right)~ r_\alpha^{-\frac{n-2k}{2}} v \left( \frac{exp_{x_\alpha}^{-1}(x)}{r_\alpha}\right) 
\end{align*}
We have that $ V_\alpha \in H_{k,0}^2({M}) $. 

\medskip\noindent{\bf Step 6:} We claim that
\begin{align}\label{cv:Va}
 V_\alpha \rightharpoonup  0 \quad \text{in} ~ H_{k,0}^2({M})\hbox{ as }\alpha\to +\infty.
\end{align}
\noindent{\it Proof of the claim}. We argue essentially as in \cite{mazumdar.gjms}. We fix $0\leq l\leq k$ and we define $\epsilon_\alpha:=r_\alpha/\tilde{r}_\alpha$ such that $\lim_{\alpha\to +\infty}\epsilon_\alpha=0$. We fix $R\geq 0$ (potentially $0$). It follows from the comparison Lemma 9.1 of \cite{mazumdar.gjms} that there exists $C>0$ such that
\begin{eqnarray*}
\int_{M\setminus B_{x_\alpha}(Rr_\alpha)}(\Delta_g^{l/2}V_\alpha)^2\, dv_g&\leq &C\int_{B_0(\delta\tilde{r}_\alpha)\setminus B_0(Rr_\alpha)}(\Delta^{l/2}(V_\alpha\circ\hbox{exp}_{x_\alpha}))^2\, dx\\
&\leq& C r_\alpha^{2(k-l)}\int_{B_0(\delta\epsilon_\alpha^{-1})\setminus B_0(R)}\left(\Delta^{l/2}\left(\eta\left(\epsilon_\alpha \cdot\right) v\right)\right)^2\, dx\\
&\leq& C r_\alpha^{2(k-l)}\int_{B_0(\delta \epsilon_\alpha^{-1})\setminus B_0(R)}|\nabla^l(\eta\left(\epsilon_\alpha \cdot\right) v)|^2\, dx\\
&\leq & C r_\alpha^{2(k-l)}\sum_{i=0}^l\int_{\rn\setminus B_0(R)}|\nabla^{l-i}\eta\left(\epsilon_\alpha \cdot\right)| |\nabla^i v|^2\, dx\\
&\leq & C r_\alpha^{2(k-l)}\sum_{i=0}^l\int_{\rn\setminus B_0(R)}\epsilon_\alpha^{2(l-i)}|\nabla^i v|^2\, dx
\end{eqnarray*}
\noindent
Since $v \in {\mathcal D}_k^{2}(\R^n)$, we have that $\nabla^iv\in {\mathcal D}_{k-i}^{2}(\R^n)$, and therefore $|\nabla^iv| \in L^{  2_{(k-i)}^{\sharp}}(\R^n)$ where $2_{(k-i)}^{\sharp}:=\frac{2n}{n-2(k-i)}$. Therefore, H\"older's inequality yields
\begin{eqnarray}
\int_{M\setminus B_{x_\alpha}(Rr_\alpha)}(\Delta_g^{l/2}V_\alpha)^2\, dv_g&\leq & C \tilde{r}_\alpha^{2(k-l)}\sum_{i=0}^l\left(\int_{\rn\setminus B_0(R)}|\nabla^i v|^{2_{(k-i)}^{\sharp}}\, dx\right)^{\frac{2}{2_{(k-i)}^{\sharp}}} \label{ineq:11}
\end{eqnarray}
Taking $R=0$ and $l=0,...,k$ yields the boundedness of $(V_\alpha)_\alpha$ in $H_{k,0}^2(M)$.

\medskip\noindent Arguing as in above, we get that for any $R>0$ and any $l=0,...,k$, we have that
\begin{equation}\label{ineq:22}
\int_{B_{x_\alpha}(Rr_\alpha)}(\Delta_g^{l/2}V_\alpha)^2\, dv_g\leq C r_\alpha^{2(k-l)}\sum_{i=0}^l\int_{B_0(R)}\epsilon_\alpha^{2(l-i)}|\nabla^i v|^2\, dx
\end{equation}
Since $\nabla^iv\in L^2_{loc}(\rn)$ for all $i=0,...,k$, then taking $l=0$ in \eqref{ineq:11} and \eqref{ineq:22}, letting $\alpha\to +\infty$ and then $R\to +\infty$ yields $V_\alpha\to 0$ in $L^2(M)$. Then the weak compactness of bounded sequences yields \eqref{cv:Va}. This proves the claim and ends Step 6.\qed

\medskip\noindent{\bf Step 7:} We claim that
\begin{equation}\label{lim:DJV}
  D J(V_\alpha) \longrightarrow 0 ~ \text{strongly as }\alpha\to +\infty
\end{equation}
\noindent{\it Proof of the claim.} We set $\varphi\in C^\infty_c(M)$. We have that 
\[
 \left\langle D J(V_\alpha), \varphi \right\rangle= 
 \int_M \Delta_g^{k/2}V_\alpha \Delta_g^{k/2}\varphi  ~{dv_g}
 - \int_{{M}} \left| V_\alpha\right|^{{2_k^{\sharp}}-2} V_\alpha \varphi ~{dv_g}
\]
We fix $R>0$ and we define 
\begin{equation*}
I_{R,\alpha}(\varphi):= \int_{B_{\xm }(Rr_\alpha)}  \Delta_g^{k/2}V_\alpha \Delta_g^{k/2}\varphi  ~{dv_g}
 - \int_{B_{\xm }(Rr_\alpha)} \left| V_\alpha\right|^{{2_k^{\sharp}}-2} V_\alpha \varphi ~{dv_g}
\end{equation*}
and
\begin{align*}
II_{R,\alpha}(\varphi):= \int_{M\setminus B_{\xm }(Rr_\alpha)}  \Delta_g^{k/2}V_\alpha \Delta_g^{k/2}\varphi  ~{dv_g}
 - \int_{M\setminus B_{\xm }(Rr_\alpha)} \left| V_\alpha\right|^{{2_k^{\sharp}}-2} V_\alpha \varphi ~{dv_g}.
\end{align*}

\medskip\noindent {\it Step 7.1:} we estimate $II_{R,\alpha}(\varphi)$. Via H\"older's and Sobolev inequality, we have that
\begin{eqnarray}\label{ineq:33} 
&&\left|II_{R,\alpha}(\varphi)\right|\leq  \left( \int_{D_\alpha(R)} (\Delta_g^{k/2}V_\alpha)^2\, {dv_g}\right)^{\frac{1}{2}}\times \Vert \Delta_g^{k/2}\varphi\Vert_{2}\\
&&+\left( \int_{ D_\alpha(R)} |V_\alpha|^{\crit}\, {dv_g}\right)^{\frac{\crit-1}{\crit}}\times \Vert \varphi\Vert_{\crit}\nonumber\\
& & \leq\left(\left( \int_{D_\alpha(R)} (\Delta_g^{k/2}V_\alpha)^2\, {dv_g}\right)^{\frac{1}{2}}+\left( \int _{D_\alpha(R)} |V_\alpha|^{\crit}\, {dv_g}\right)^{\frac{\crit-1}{\crit}}\right)\cdot \Vert \varphi\Vert_{H_k^2}\nonumber
\end{eqnarray}
with $D_\alpha(R):=M\setminus B_{x_\alpha}(R r_\alpha)$. Lemma 9.1 in \cite{mazumdar.gjms} and $v\in L^{\crit}(\rn)$ yield 
\begin{equation}\label{bnd:Vm}
\int_{M\setminus B_{x_\alpha}(R r_\alpha)} |V_\alpha|^{\crit}\, {dv_g}\leq C \int_{\rn\setminus B_{0}(R r_\alpha) } |V_\alpha\circ\exp_{x_\alpha}|^{\crit}\, dx\leq C\int_{\rn\setminus B_{0 }(R)} |v|^{\crit}\, dx
\end{equation}
Plugging \eqref{ineq:11} with $l=k$ and \eqref{bnd:Vm} into \eqref{ineq:33}, letting $R\to +\infty$ and $\alpha\to +\infty$ yields
\begin{equation}\label{eq:55}
\lim_{R\to +\infty}\lim_{\alpha\to +\infty}\frac{II_{R,\alpha}(\varphi)}{\Vert\varphi\Vert_{H_k^2}}=0\hbox{ uniformly wrt }\varphi\in H_{k,0}^2(M)\setminus\{0\}
\end{equation} 

\medskip\noindent{\bf Step 7.2:} We now estimate $I_{R,\alpha}(\varphi)$. We define
  \begin{equation*}
   \overline{\varphi}_\alpha(x) = \eta(\epsilon_\alpha x) r_\alpha^{\frac{n-2k}{2}} \varphi \left( exp_{x_\alpha}(r_\alpha x) \right)
  \end{equation*}
where $\epsilon_\alpha:=r_\alpha/\tilde{r}_\alpha$. As one checks, $\overline{\varphi}_\alpha\in C^\infty_c(\rn)$. Using the comparison Lemma 9.1 in \cite{mazumdar.gjms} and arguing as in \eqref{ineq:11}-\eqref{ineq:22}, we get that
\begin{equation*}
\Vert \overline{\varphi}_\alpha\Vert_{{\mathcal D}_k^2(\rn)}\leq C \Vert \varphi\Vert_{H_k^2}
\end{equation*}
where $C>0$ is independent of $\varphi$. As one checks,
\begin{equation*}
I_{R,\alpha}(\varphi)= \int_{B_0(R) } 
  \Delta_{{g}_\alpha}^{k/2}v \Delta_{{g}_\alpha}^{k/2} \overline{\varphi}_\alpha \, {dv_{{g}_\alpha}}-\int_{B_0(R) } 
  |v|^{\crit-2}v \overline{\varphi}_\alpha \,{dv_{{g}_\alpha}}
\end{equation*}
Since $g_\alpha\to \hbox{Eucl}$ as $\alpha\to +\infty$ in $C^p_{loc}(\rn)$ for all $p\geq 1$, we get 
\begin{equation}\label{ineq:44}
I_{R,\alpha}(\varphi)= \int_{B_0(R) } 
  \Delta^{k/2}v \Delta^{k/2} \overline{\varphi}_\alpha \, dx-\int_{B_0(R) } 
  |v|^{\crit-2}v \overline{\varphi}_\alpha \,dx+ o\left(\Vert \overline{\varphi}_\alpha\Vert_{{\mathcal D}_k^2(\rn)}\right)
  \end{equation}
where the convergence is uniform wrt $\overline{\varphi}_\alpha$. Since $v$ is a weak solution to \eqref{eqn1}, then \eqref{ineq:44} yields
\begin{equation}\label{eq:66}
\lim_{R\to +\infty}\lim_{\alpha\to +\infty}\frac{I_{R,\alpha}(\varphi)}{\Vert\varphi\Vert_{H_k^2}}=0\hbox{ uniformly wrt }\varphi\in H_{k,0}^2(M)\setminus\{0\}
\end{equation} 

\medskip\noindent The limits \eqref{eq:55} and \eqref{eq:66} yield $\langle DJ( V_\alpha),\varphi\rangle=o(\Vert\varphi\Vert_{H_k^2})$ as $\alpha\to +\infty$ uniformly wrt $\varphi\in C^\infty_c(M)$. The boundedness of $(V_\alpha)$ in $H_{k,0}^2(M)$ then yields $DJ( V_\alpha)\to 0$ strongly in $(H_{k,0}^2(M))^\prime$ as $\alpha\to +\infty$. This proves \eqref{lim:DJV} and ends Step 7.\qed

\medskip\noindent We define $w_\alpha:= v_\alpha- V_\alpha$. It follows from \eqref{cv:Va} that $w_\alpha \rightharpoonup 0$ weakly in $H_{k,0}^{2}(M)$. 

\medskip\noindent{\bf Step 8:} We claim that
\begin{equation}\label{79}
  D J(w_\alpha) \longrightarrow 0 ~ \text{strongly}
\end{equation}
{\it Proof of the claim.} For $\varphi \in H_{k,0}^2(M)$, we write
\begin{equation}\label{80}
 \left\langle DJ(w_\alpha),\varphi\right\rangle =
 \left\langle DJ(v_\alpha),\varphi\right\rangle - \left\langle DJ(V_\alpha),\varphi\right\rangle -
 \int_M  \varPhi_\alpha \varphi ~{dv_g}  
\end{equation}
where $\varPhi_\alpha :=  \left| w_\alpha\right|^{ 2_k^{\sharp}-2} w_\alpha     -  \left| v_\alpha \right|^{ 2_k^{\sharp}-2} v_\alpha   +
  \left| V_\alpha \right|^{ 2_k^{\sharp}-2}V_\alpha $. Then by applying the  H\"{o}lder and Sobolev inequalities  we get
\begin{equation*}
 \left| ~ \int_M  \varPhi_\alpha \varphi ~ {dv_g} ~  \right| \leq  C \left\| \varphi \right\|_{H_k^2}
 \left\| \varPhi_\alpha   \right\|_{{{2_k^{\sharp}}}/{({2_k^{\sharp} -1)}}}
\end{equation*}

\smallskip\noindent{\it Step 8.1:} We fix $R>0$. Inequality \eqref{ineq:useful} and H\"older's inequality yield
\begin{eqnarray*}
&&\int_{M\setminus B_{x_\alpha}(R r_\alpha)}|\varPhi_\alpha|^{\crit/(\crit-1)}\, dv_g\\
&&\leq C\int_{M\setminus B_{x_\alpha}(R r_\alpha)}\left(|v_\alpha|^{\crit-2}|V_\alpha|+|V_\alpha|^{\crit-2}|v_\alpha|\right)^{\crit/(\crit-1)}\, dv_g\\
 &&\leq C\left(\int_M |v_\alpha|^{\crit}\, dv_g\right)^{\frac{\crit-2}{\crit-1}} \left(\int_{M\setminus B_{x_\alpha}(R r_\alpha)} |V_\alpha|^{\crit}\, dv_g\right)^{\frac{1}{\crit-1}}\\
&& + C\left(\int_M |v_\alpha|^{\crit}\, dv_g\right)^{\frac{1}{\crit-1}} \left(\int_{M\setminus B_{x_\alpha}(R r_\alpha)} |V_\alpha|^{\crit}\, dv_g\right)^{\frac{\crit-2}{\crit-1}}
\end{eqnarray*}
Since $(v_\alpha)$ is uniformly bounded in $H_k^2(M)$, then \eqref{bnd:Vm} yields
\begin{equation}\label{85}
\lim_{R\to +\infty}\lim_{\alpha\to +\infty}\int_{M\setminus B_{x_\alpha}(R r_\alpha)}|\varPhi_\alpha|^{\crit/(\crit-1)}\, dv_g=0.
\end{equation}
This ends Step 8.1.

\smallskip\noindent{\it Step 8.2:} We fix $R>0$. A change of variable and inequality \eqref{ineq:useful} yield
\begin{eqnarray*}
&&\int_{B_{x_\alpha}(R r_\alpha)}|\varPhi_\alpha|^{\crit/(\crit-1)}\, dv_g\\
&&= \int_{B_0(R)}\left|\left| \tilde{v}_\alpha-v\right|^{ 2_k^{\sharp}-2} (\tilde{v}_\alpha-v)     -  \left| \tilde{v}_\alpha \right|^{ 2_k^{\sharp}-2} \tilde{v}_\alpha   +
  \left| v \right|^{ 2_k^{\sharp}-2}v\right|^{\crit/(\crit-1)}\, dv_{g_\alpha}\\
  &&\leq  C\int_{B_0(R)}\left(\left| \tilde{v}_\alpha-v\right|^{ \frac{(\crit-2)\crit}{\crit-1}} |v|^{ \frac{\crit}{\crit-1}}  + \left| v\right|^{ \frac{(\crit-2)\crit}{\crit-1}}\left| \tilde{v}_\alpha-v\right|^{ \frac{\crit}{\crit-1}} \right)\, dx
\end{eqnarray*}
For any $\eta\in C^\infty_c(\rn)$, we have that $\eta\tilde{v}_\alpha\rightharpoonup \eta v$ weakly in ${\mathcal D}_k^2(\rn)$. Therefore, up to extracting a subsequence, $(\tilde{v}_\alpha)_\alpha$ is uniformly bounded in $L^{\crit}(B_0(R))$ and goes to $v$ almost everywhere as $\alpha\to +\infty$. Therefore Lemma \ref{int.lemma} yields that for any $R>0$,
\begin{equation}\label{86}
\lim_{\alpha\to +\infty}\int_{B_{x_\alpha}(R r_\alpha)}|\varPhi_\alpha|^{\crit/(\crit-1)}\, dv_g=0.
\end{equation} 

\medskip\noindent The limits \eqref{85}-\eqref{86} yield  $ \left\| \varPhi_\alpha   \right\|_{{{2_k^{\sharp}}}/{({2_k^{\sharp} -1)}}}\to 0$ as $\alpha\to +\infty$. Then  by \eqref{80}  we get $DJ(w_\alpha)\to 0$ in $(H_{k,0}^2(M))^\prime$ as $\alpha\to +\infty$. This proves \eqref{79} and ends Step 8.\qed

\medskip\noindent{\bf Step 9:} We claim that we have the following decomposition of energy.
 \begin{equation}\label{est:nrj}
  J(w_\alpha)= J(v_\alpha)- E(v)+ o(1)  \hbox{ where }o(1) \rightarrow 0\hbox{ as }\alpha \rightarrow + \infty.
 \end{equation}

 \medskip\noindent{\it Proof of the claim}. As one checks, 
 \begin{eqnarray*}
&&J(v_\alpha)-J(w_\alpha)- J(V_\alpha)=\langle DJ(w_\alpha),V_\alpha\rangle\\
&&-\frac{1}{\crit}\int_M\left(|w_\alpha+V_\alpha|^{\crit}-|w_\alpha|^{\crit}-\crit|w_\alpha|^{\crit-2}w_\alpha V_\alpha-|V_\alpha|^{\crit}\right) \, dv_g \end{eqnarray*}
We fix $R>0$. Arguing as in the proof of \eqref{86}, we get that
\begin{equation*}
\lim_{\alpha\to +\infty}\int_{B_{x_\alpha}(R r_\alpha)}\left(|w_\alpha+V_\alpha|^{\crit}-|w_\alpha|^{\crit}-\crit|w_\alpha|^{\crit-2}w_\alpha V_\alpha-|V_\alpha|^{\crit}\right) \, dv_g =0.
\end{equation*}
As one checks, there exists $C>0$ such that 
$$\left| |a+b|^{\crit}-|a|^{\crit}-\crit |a|^{\crit-2}ab-|b|^{\crit}\right|\leq C\left(|a|^{\crit-2}|b|^2+|a|\cdot |b|^{\crit-1}\right)$$
 for all $a,b\in\R$. As in the proof of \eqref{85}, we get that 
\begin{equation*}
\lim_{R\to +\infty}\lim_{\alpha\to +\infty}\int_{D_\alpha(R)}\left(|w_\alpha+V_\alpha|^{\crit}-|w_\alpha|^{\crit}-\crit|w_\alpha|^{\crit-2}w_\alpha V_\alpha-|V_\alpha|^{\crit}\right) \, dv_g =0, 
\end{equation*}
where $D_\alpha(R):=M\setminus B_{x_\alpha}(R r_\alpha)$. These yield $ J(v_\alpha)=J(w_\alpha)+ J(V_\alpha)+o(1)$.

\medskip\noindent We now estimate $J(V_\alpha)$. The estimates \eqref{ineq:11} and \eqref{bnd:Vm} yield
\begin{equation*}
\lim_{R\to +\infty}\lim_{\alpha\to +\infty}\int_{M\setminus B_{x_\alpha}(R r_\alpha)}\left((\Delta_g^{k/2} V_\alpha)^2+|V_\alpha|^{\crit}\right)\, dv_g=0
\end{equation*}
For $R>0$, we have that
\begin{equation*}
\int_{B_{x_\alpha}(R r_\alpha)}\left(\frac{(\Delta_g^{k/2} V_\alpha)^2}{2}-\frac{|V_\alpha|^{\crit}}{\crit}\right)\, dv_g=\int_{B_{0}(R)}\left(\frac{(\Delta_{g_\alpha}^{k/2} v)^2}{2}-\frac{|v|^{\crit}}{\crit}\right)\, dv_{g_\alpha}
\end{equation*}
Since $g_\alpha\to\hbox{Eucl}$ locally uniformly in $C^p$ for all $p$ and $v\in {\mathcal D}_k^2(\rn)$, we get that
\begin{equation*}
\lim_{R\to +\infty}\lim_{\alpha\to +\infty} \int_{B_{x_\alpha}(R r_\alpha)}\left(\frac{(\Delta_g^{k/2} V_\alpha)^2}{2}-\frac{|V_\alpha|^{\crit}}{\crit}\right)\, dv_g=\int_{\rn}\left(\frac{(\Delta^{k/2} v)^2}{2}-\frac{|v|^{\crit}}{\crit}\right)\, dx
\end{equation*}
All these estimates yield \eqref{est:nrj}. This ends Step 9.\qed

\medskip\noindent{\bf Step 10:} Next we deal with the case
 \begin{equation*}
 d_g (x_\alpha, \partial M)  =O(r_\alpha)~ \text{ as  ${\alpha \rightarrow + \infty}$}
  \end{equation*}
Since $r_\alpha\to 0$ as $\alpha\to +\infty$, then there exists $x_\infty\in\partial M$ such that $x_\alpha \rightarrow x_{\infty}$ as $\alpha\to +\infty$. For any $\alpha\in\N$, we let $z_{\alpha} \in \partial M$ be such that 
\begin{align*}
d_{g}(x_{\alpha}, z_{\alpha})= d_g (x_\alpha, \partial M)
\end{align*}
In particular, $\lim_{\alpha\to +\infty}z_\alpha=x_\infty$. We choose a family of charts $z\mapsto {\mathcal T}_z$ for $z\in \Omega\cap\partial M$ as in \eqref{def:T}. Since the $d({\mathcal T}_z)_0$ is an isometry, there exists  ${\mathcal{C}}_{1}, {\mathcal{C}}_{2}>0$, $\tau_{1}, \tau_{2}>0$ such that for any $z\in \Omega \cap\partial M$, $r < \tau_{1}$ and $y \in \R^{n}_{-}\cap B_0(\tau_2)$, one has
\begin{align*}
B_{{\mathcal{T}_{z}}(y)}({\mathcal{C}}_{1} r) \cap M  \subset {\mathcal{T}_{z}} \left( B_{y} (r) \cap \R_{-}^{n}  \right) \subset  B_{{\mathcal{T}_{z}}(y)}({\mathcal{C}}_{2} r) \cap M 
\end{align*}
For $ x \in  r_\alpha^{-1} U  \cap \{x_{1} <0\} $, we define 
$$\tilde{v}_{\alpha}(x):=  r_{\alpha}^{\frac{n-2k}{2}} v_\alpha \circ \mathcal{T}_{z_\alpha}(r_{\alpha}x ) \hbox{ and }\tilde{g}_{\alpha}(x) := {\mathcal{T}_{z_\alpha}}^{\star} g~ (r_{\alpha}x  )$$
As one checks, for any $\eta\in C^\infty_c(\rn)$, we have that $\eta \tilde{v}_\alpha\in {\mathcal D}_k^2(\rnm)$. Arguing as Step 4, we get that there exists $v\in {\mathcal D}_k^2(\rnm)$ such that
\begin{equation*}
\eta  \tilde{v}_\alpha\rightharpoonup \eta v\hbox{ weakly in }{\mathcal D}_k^2(\rnm)\hbox{ as }\alpha\to +\infty.
\end{equation*}
Moreover, using Proposition \ref{prop:fund} and arguing as in Step 5, we get that $v\not\equiv 0$ is a weak solution to \eqref{lim:2} and $\tilde{v}_\alpha\to v$ as $\alpha\to +\infty$ strongly in $H_{k}^2(B_0(R)\cap \rnm)$ for all $R>0$. As in Definition \ref{def:bubbles}, for $\alpha\in\N$ and $x\in \overline{M}$, we set 
 \begin{equation*}
 V_\alpha(x):=  B_{z_\alpha,r_\alpha}(v) (x)=
\eta\left(\mathcal{T}_{z_\alpha}^{-1}(x)\right)   r_\alpha^{-\frac{n-2k}{2}} v \left(r_\alpha^{-1} {\mathcal{T}_{z_\alpha}}^{-1}(x) \right)
 \end{equation*} 
We define $w_{\alpha}:= v_{\alpha} - V_{\alpha}$. Arguing as in Steps 6 to 9, we get that

\begin{enumerate}
\item[$\bullet$] $w_\alpha \rightharpoonup 0 \hbox{ weakly in } H_{k,0}^{2}(M)$ 
\item[$\bullet$] $ D J(w_\alpha) \to 0  \hbox{ weakly in } (H_{k,0}^{2}(M))^\prime$
\item[$\bullet$] $ J(w_\alpha)= J(v_\alpha)- E(v)+ o(1) $
\end{enumerate}
as $\alpha\to +\infty$. This completes  the proof of Lemma \ref{bubble.lemma}.

\section{Nonnegative Palais-Smale sequences}\label{sec:positive}
To prove Theorem \ref{th:2}, we first set the following property:

\begin{proposition}\label{prop:hr}  Let $(\ua)$ be  a Palais-Smale sequence for the functional $I$ on the space $ H_{k,0}^2(M)$. Let $d\in\N$ and $d$ bubbles $[(x_\alpha^{(j)}), (r_\alpha^{(j)}), u^{(j)}]$, $j=1,...,d$, be as in Theorem \ref{th:1}. Then, for any $N \in \{1,\ldots ,d \}$, there exists $L\geq 0$ sequences $(y_{\alpha}^{j})_{\alpha >0}\in \overline{M}$ and $(\lambda_{\alpha}^{j})_{\alpha >0}\in (0,+\infty)$, $j=1, \cdots, L$, such that for any $R >0$ 
\begin{equation*}
\lim_{R'\to +\infty}\lim_{\alpha\to +\infty}\int_{\left( B_{x_{\alpha}^{N}}(R r_{\alpha}^{N}) \backslash \bigcup_{j=1}^{L} B_{y_{\alpha}^{j}}(R^{\prime} \lambda_{\alpha}^{j}) \right) \cap M}  | \ua-B_{x_{\alpha}^{(N)}, r_{\alpha}^{(N)}} ( u^{(N)})|^{2^{\sharp}_{k}}~dv_{g} =0
\end{equation*}
where for any  $j$, $j=1, \cdots, L$, $d_{g}(x_{\alpha}^{N}, y_{\alpha}^{j})= o(r_{\alpha}^{N})$ and $\lambda_{\alpha}^{j}=o(r_{\alpha}^{N})$ as $\alpha\to +\infty$. Moreover, we have that
$$\lim_{\alpha\to +\infty}\frac{d_g(x_\alpha^i,x_\alpha^j)^2}{r_\alpha^i r_\alpha^j}+\frac{r_\alpha^i}{r_\alpha^j}+\frac{r_\alpha^j}{r_\alpha^i}=+\infty\hbox{ for all }i\neq j\in \{1,...,d\}.$$
\end{proposition}
We omit the proof that goes exactly as in Hebey-Robert \cite{Hebey-robert}, by using the boundary chart \eqref{def:T} for bubbles accumulating on the boundary.

\medskip\noindent We now prove Theorem \ref{th:2}. We let $(\ua)_\alpha$ be as in the statement of the theorem, and we let $[(x_\alpha^{(j)}), (r_\alpha^{(j)}), u^{(j)}]$, $j=1,...,d$, be the associated bubbles. We fix $N\in \{1,...,d\}$. For simplicity, we define $r_\alpha:=r_\alpha^{(N)}$ and $x_\alpha:=x_\alpha^{(N)}$. We assume that $r_\alpha^{-1}d(x_\alpha,\partial M)\to +\infty$ as $\alpha\to +\infty$. It then follows from Proposition \ref{prop:hr}  that there exists a finite set $\mathcal S\subset \rn$ such that $\lim_{\alpha\to +\infty}\tva=u^N$ strongly in $L^{\crit}_{loc}(\rn\setminus {\mathcal S})$ where $\tva(x):=r_\alpha^{\frac{n-2k}{2}}\ua(\hbox{exp}_{x_\alpha}(r_\alpha x))$ for $x\in \rn$. Up to extracting, the convergence holds a.e. Since $\ua\geq 0$, we then get that $u^N\geq 0$. It then follows from Lemma 4 in Ge-Wei-Zhou \cite{poly} that there exists $\lambda>0$ and $a\in\rn$ such that $u^N=U_{\lambda,a}$ is of the form \eqref{eq:ext}. 

\medskip\noindent We claim that  $u^N=U_{\lambda,0}$, that is $a=0$. We prove the claim. Indeed, rescaling \eqref{def:mu} and \eqref{def:ra:xa} yields
$$\int_{r_\alpha^{-1}\hbox{exp}_{x_\alpha}^{-1}(B_{\hbox{exp}_{x_\alpha}(r_\alpha x)}(r_\alpha))}|\tva|^{\crit}\, dv_{g_\alpha}\leq \int_{B_0(1)}|\tva|^{\crit}\, dv_{g_\alpha}$$
for all $z\in\rn$ and $\alpha$ large enough. Since the exponential is a normal chart and isometric at $x_\alpha$, we get that for all $z\in\rn$ and all $\epsilon>0$
$$\hbox{exp}_{x_\alpha}\left(r_\alpha B_z(1-\epsilon)\right)\subset B_{\hbox{exp}_{x_\alpha}(r_\alpha z)}(r_\alpha).$$
Plugging these two inequalities together, letting $\alpha\to +\infty$, using the strong convergence \eqref{strong:cv}, we get that $\int_{B_z(1-\epsilon)}|u^N|^{\crit}\, dx\leq \int_{B_0(1)}|u^N|^{\crit}\, dx$. Letting $\epsilon\to 0$ yields
$$\int_{B_z(1)}|u^N|^{\crit}\, dx\leq \int_{B_0(1)}|u^N|^{\crit}\, dx.$$
As one checks, since $u^N=U_{\lambda,a}$ is as in \eqref{eq:ext}, the maximum of the left-hand-side is achieved if and only if $z=a$. Therefore $a=0$ and $u^N=U_{\lambda,0}$. This proves the claim.

\smallskip\noindent As a consequence, as one checks, when $r_\alpha^{-1}d(x_\alpha,\partial M)\to +\infty$ as $\alpha\to +\infty$, the bubble rewrites
$$B_{x_\alpha,r_\alpha}(u^N)=B_{x_\alpha, \lambda r_\alpha}(U_{1,0})=\eta\left(\frac{\hbox{exp}_{x_\alpha}^{-1}(\cdot)}{\tilde{r}_\alpha}\right)\alpha_{n,k}\left(\frac{\lambda r_\alpha}{\lambda^2r_\alpha^2+d_g(\cdot,x_\alpha)^2}\right)^{\frac{n-2k}{2}}.$$

\smallskip\noindent We fix $N\in \{1,...,d\}$. We claim that $(r_\alpha^N)^{-1}d(x_\alpha^N,\partial M)\to +\infty$ as $\alpha\to +\infty$. We argue by contradiction and we assume that the limit is finite. We argue as in the case above. Up to rescaling, and using the boundary chart \eqref{def:T}, we get that $\ua$ goes to $u^N$ strongly as $\alpha\to +\infty$ in $L^{\crit}_{loc}(\rn\setminus {\mathcal S})$, where $\mathcal S$ is finite. Therefore $u^N$ is a nonegative nonzero weak solution to \eqref{lim:2}, contradicting Lemma 3 in Ge-Wei-Zhou \cite{poly}. Therefore the limit is infinite and we are back to the previous case. 

\smallskip\noindent All these steps prove Theorem \ref{th:2}.

\end{document}